\documentstyle[twoside]{article}
\title{Exponentially small expansions related to the parabolic cylinder function}
\author{\sc R. B.\ Paris \\
{\em Division of Computing and Mathematics,}\\
{\em Abertay University, Dundee DD1 1HG, UK}}
\begin{document}
\def\f#1#2{\mbox{${\textstyle \frac{#1}{#2}}$}}
\def\dfrac#1#2{\displaystyle{\frac{#1}{#2}}}
\newcommand{\bee}{\begin{equation}}
\newcommand{\ee}{\end{equation}}
\newcommand{\fr}{\frac{1}{2}}
\newcommand{\fs}{\f{1}{2}}
\newcommand{\ff}{\f{1}{4}}
\newcommand{\g}{\Gamma}
\newcommand{\br}{\biggr}
\newcommand{\bl}{\biggl}
\newcommand{\ra}{\rightarrow}
\newcommand{\gtwid}{\raisebox{-.8ex}{\mbox{$\stackrel{\textstyle >}{\sim}$}}}
\newcommand{\ltwid}{\raisebox{-.8ex}{\mbox{$\stackrel{\textstyle <}{\sim}$}}}
\renewcommand{\topfraction}{0.9}
\renewcommand{\bottomfraction}{0.9}
\renewcommand{\textfraction}{0.05}
\newcommand{\mcol}{\multicolumn}
\date{}
\maketitle
\pagestyle{myheadings}
\markboth{\hfill \sc R. B.\ Paris  \hfill}
{\hfill \sc Exponentially small asymptotic expansion \hfill}
\begin{abstract}
The refined asymptotic expansion of the confluent hypergeometric function $M(a,b,z)$ on the Stokes line $\arg\,z=\pi$ given in {\it Appl. Math. Sci.} {\bf 7} (2013) 6601--6609 is employed to derive the correct exponentially small contribution to the asymptotic expansion for the even and odd solutions of a second-order differential equation related to Weber's equation. It is demonstrated that the standard asymptotics of the parabolic cylinder function $U(a,z)$ yield an incorrect exponentially small contribution to these solutions. Numerical results verifying the accuracy of the new expansions are given.  
\vspace{0.4cm}

\vspace{0.3cm}

\end{abstract}
\vspace{0.3cm}

\noindent $\,$\hrulefill $\,$

\vspace{0.2cm}

\begin{center}
{\bf 1. \  Introduction}
\end{center}
\setcounter{section}{1}
\setcounter{equation}{0}
\renewcommand{\theequation}{\arabic{section}.\arabic{equation}}
The second-order differential equation
\bee\label{e11}
\frac{d^2w}{dx^2}+x\,\frac{dw}{dx}+(a+\fs)w=0,\qquad a\in{\bf R}
\ee
is a transformation of Weber's equation. The solutions of (\ref{e11}) are $e^{-x^2/4} U(a,\pm ix)$, where $U(a,z)$ is the parabolic cylinder function defined in terms of the confluent hypergeometric functions of the first ($M(a,b,z)$) and second ($U(a,b,z)$) kinds by \cite[(12.7.16)]{DLMF}
\begin{eqnarray*}
U(a,z)&=&2^{-\frac{1}{4}-\frac{1}{2}a} e^{-z^2/4} U(\fs a+\f{1}{4}, \fs, \fs z^2)\\
&=&2^{-\frac{1}{4}-\frac{1}{2}a}\sqrt{\pi}e^{-z^2/4}\bl\{\frac{M(\fs a+\f{1}{4},\fs,\fs z^2)}{\g(\fs a+\f{3}{4})}-2^\fr z\ \frac{M(\fs a+\f{3}{4},\f{3}{2},\fs z^2)}{\g(\fs a+\f{1}{4})}\br\}.
\end{eqnarray*}
Even and odd solutions of (\ref{e11}) satisfying  the initial conditions $w(0)=1$, $w'(0)=0$ and $w(0)=0$, $w'(0)=1$ are respectively
\bee\label{e12}
w_1(a,x)=M(\fs a+\ff,\fs,-\fs x^2),\qquad w_2(a,x)=x M(\fs a+\f{3}{4}, \f{3}{2},-\fs x^2),
\ee
%where $M(a,b,z)$ is the confluent hypergeometric function of the first kind (or first Kummer function,
%also denoted by ${}_1F_1(a; b; z)$). 
From \cite[(13.6.14,15)]{DLMF}, these solutions can be expressed alternatively in terms of the parabolic cylinder functions $U(a,\pm ix)$ by
\bee\label{e13}
\left.
\begin{array}{lll}
w_1(a,x)&=&\dfrac{2^{\frac{1}{2}a-\frac{3}{4}}}{\sqrt{\pi}} \g(\fs a+\f{3}{4}) e^{-x^2/4} \{U(a,ix)+U(a,-ix)\}\\
\\
w_2(a,x)&=&\dfrac{2^{\frac{1}{2}a-\frac{5}{4}}}{\sqrt{\pi}}\g(\fs a+\f{1}{4}) e^{-x^2/4} \{U(a,-ix)-U(a,ix)\}.
\end{array}\right\}
\ee

For half-integer values of the parameter $a$ the solutions $w_1(a,x)$ and $w_2(a,x)$ reduce to simpler forms.
For non-negative integer $n$, the solutions $w_1(a,x)$ and $w_2(a,x)$ reduce to a polynomial multiplied by $\exp\,[-x^2/2]$ when $a=2n+\fs$ and $a=2n+\f{3}{2}$ respectively, whereas they reduce to polynomials when $a=-2n-\fs$ and $a=-2n-\f{3}{2}$. Apart from these special values our interest in this note will be concerned with the asymptotic nature of the solutions in (\ref{e13}) as $x\to+\infty$ for real values of $a$, paying particular attention to the exponentially small contributions.

The interest in exponentially precise asymptotics  during the past three decades has shown that retention of exponentially small terms, previously neglected in asymptotics, is essential for a high-precision description. An early example  that illustrated the advantage of retaining exponentially small terms in the asymptotic expansion of a certain integral was given in Olver's well-known book \cite[p.~76]{Ob}.
Although such subdominant terms are negligible in the Poincar\'e sense, their inclusion can significantly improve the numerical accuracy in computations.

Use of the standard asymptotic expansion of $U(a,\pm ix)$ for $x\to+\infty$ and fixed $a$ shows that $w_1(a,x)$ and $w_2(a,x)$ possess an algebraic-type expansion of $O(x^{-a-1/2})$, together with an exponentially small contribution of $O(x^{a-1/2} \exp\,[-x^2/2])$. However, it will be shown that, when the algebraic expansion is optimally truncated at, or near, its least term (in magnitude), the exponentially small contribution obtained in this manner is out by a factor 2. And furthermore it will be shown that there is an additional expansion present in the exponentially small contribution that is not accounted for in the standard asymptotics. These conclusions follow from an application to the representations in (\ref{e12}) of the recent treatment of the asymptotic expansion of the Kummer 
function $M(a,b,z)$ presented in \cite{P13}, which takes account of the Stokes phenomenon on the negative $z$-axis.

\vspace{0.6cm}

\begin{center}
{\bf 2. \ The standard asymptotic expansion for $w_r(a,x)$ ($r=1, 2$)}
\end{center}
\setcounter{section}{2}
\setcounter{equation}{0}
\renewcommand{\theequation}{\arabic{section}.\arabic{equation}}
We first derive the asymptotic expansions of the solutions $w_1(a,x)$ and $w_2(a,x)$ for $x\to+\infty$ and fixed $a$ employing the standard large-$z$ expansion of $U(a,z)$. It will be demonstrated in Section 3 that the coefficient multiplying the exponentially small contribution obtained in this manner is out by a factor of 2.

From \cite[(12.9.3)]{DLMF}, the expansion of $e^{z^2/4} U(a,z)$ as $|z|\to\infty$ is given by
\[e^{z^2/4} U(a,z) \hspace{10cm}\]
\bee\label{e21}
\sim z^{-a-\fr} \sum_{k=0}^\infty \frac{(-)^k (\fs+a)_{2k}}{k! (2z^2)^k}\pm \frac{i\sqrt{2\pi}\,e^{\mp\pi ia}}{\g(\fs+a)}\, z^{a-\fr} e^{-z^2/2} \sum_{k=0}^\infty \frac{(\fs-a)_{2k}}{k! (2z^2)^k}
\ee
\[\hspace{8cm}(\f{1}{4}\pi<\pm\arg\,z<\f{5}{4}\pi),\]
where $(\alpha)_k=\g(\alpha+k)/\g(\alpha)$ is the Pochhammer symbol.
Then some routine algebra shows that
\[e^{-x^2/4}\{U(a,ix)+U(a,-ix)\}\hspace{8cm}\]
\[\sim
2\sin \pi\vartheta\, x^{-a-\fr} \sum_{k=0}^\infty \frac{(\fs+a)_{2k}}{k! (2x^2)^k}
+\frac{2\sqrt{2\pi}}{\g(\fs+a)}\,\cos \pi\vartheta\,x^{a-\fr} e^{-x^2/2} \sum_{k=0}^\infty \frac{(-)^k(\fs-a)_{2k}}{k! (2x^2)^k}\] 
and
\[e^{-x^2/4}\{U(a,-ix)-U(a,ix)\}\hspace{8cm}\]
\[\sim
2i\cos \pi\vartheta\, x^{-a-\fr} \sum_{k=0}^\infty \frac{(\fs+a)_{2k}}{k! (2x^2)^k}
+\frac{2i\sqrt{2\pi}}{\g(\fs+a)}\,\sin \pi\vartheta\,x^{a-\fr} e^{-x^2/2} \sum_{k=0}^\infty \frac{(-)^k(\fs-a)_{2k}}{k! (2x^2)^k}\]
as $x\to+\infty$, where we have defined
\bee\label{e22}
\vartheta:=\frac{1}{2} a-\frac{1}{4}.
\ee

From (\ref{e13}) followed by use of the duplication formula for the gamma function, we then obtain
the expansions
\[w_1(a,x)\sim \frac{2^{\frac{1}{4}+\fr a}\sqrt{\pi}}{\g(\f{1}{4}-\fs a)}\,x^{-a-\fr} \sum_{k=0}^\infty \frac{(\fs+a)_{2k}}{k! (2x^2)^k}\hspace{4cm}\]
\bee\label{e24a}
\hspace{3cm}+\frac{2^{\frac{1}{4}-\fr a}\sqrt{\pi}}{\g(\f{1}{4}+\fs a)}
\,2\cos \pi\vartheta\,x^{a-\fr} e^{-x^2/2} \sum_{k=0}^\infty \frac{(-)^k(\fs-a)_{2k}}{k! (2x^2)^k}\ee
and
\[w_2(a,x)\sim \frac{2^{-\frac{1}{4}+\fr a}\sqrt{\pi}}{\g(\f{3}{4}-\fs a)}\,x^{-a-\fr} \sum_{k=0}^\infty \frac{(\fs+a)_{2k}}{k! (2x^2)^k}\hspace{4cm}\]
\bee\label{e24b}
\hspace{3cm}+\frac{2^{-\frac{1}{4}-\fr a}\sqrt{\pi}}{\g(\f{3}{4}+\fs a)}
\,2\sin \pi\vartheta\,x^{a-\fr} e^{-x^2/2} \sum_{k=0}^\infty \frac{(-)^k(\fs-a)_{2k}}{k! (2x^2)^k}\ee
for $x\to+\infty$. When $a=\fs$, the expression in (\ref{e12}) shows that $w_1(a,x)=\exp\,[-x^2/2]$, since the confluent hypergeometric function reduces to an exponential in this case. But (\ref{e24a}) yields $w_1(\fs,x)\sim 2 \exp\,[-x^2/2]$. Similarly, when $s=\f{3}{2}$ we see from (\ref{e12}) that $w_2(\f{3}{2},x)=x \exp\,[-x^2/2]$ but (\ref{e24b}) yields $w_2(\f{3}{2},x)\sim 2x \exp\,[-x^2/2]$. In both cases
the standard asymptotic expansions predict an exponentially small term that is out by a factor of 2.

In the next section we shall show how the use of the more refined asymptotics of $M(a,b,z)$ on the negative $z$-axis yields a more accurate exponentially small contribution. 

\vspace{0.6cm}

\begin{center}
{\bf 3. \ The more refined expansion for $w_r(a,x)$ ($r=1, 2$) as $x\ra+\infty$}
\end{center}
\setcounter{section}{3}
\setcounter{equation}{0}
\renewcommand{\theequation}{\arabic{section}.\arabic{equation}}
In \cite{P13}, the expansion of the Kummer function $M(a,b,-x)$ for $x\to+\infty$ 
was re-considered. This treatment took into account the Stokes phenomenon on the negative real axis to yield the correct asymptotic behaviour of the exponentially small contribution. The resulting expansion has the form \cite[Theorem 1]{P13}:
\[\frac{\g(a)}{\g(b)}\,M(a,b,-x)-\frac{x^{-a}\g(a)}{\g(b-a)} \sum_{k=0}^{m_0-1}\frac{(a)_k(1+a-b)_k}{k! x^k} \hspace{4cm}\]
\bee\label{e31}
=x^{\vartheta} e^{-x}\bl\{ \cos \pi{\hat\vartheta} \sum_{j=0}^{M-1}(-)^jA_j\,x^{-j}-\frac{2\sin \pi{\hat\vartheta}}{\sqrt{2\pi x}}\sum_{j=0}^{M-1}(-)^jB_j\,x^{-j}+O(x^{-M})\br\}
\ee
as $x\to+\infty$, where ${\hat\vartheta}=a-b$, $M$ denotes a positive integer and for simplicity we shall suppose that $a$, $b$ are real parameters\footnote{The expansion (\ref{e31}) also holds when $\vartheta$ is a negative integer; see \cite[Theorem 2]{P13}. When $\vartheta=n$, a positive integer, the algebraic expansion vanishes to leave just the first series in the exponentially small contribution with the coefficients $A_j=0$ for $j>n$.} 
Here, the dominant algebraic expansion on the left-hand side of (\ref{e31}) has been optimally truncated with index $m_o$ given by
\bee\label{e33}
m_0=x-2a+b+\alpha,\qquad |\alpha|<1.
\ee

The expansion on the right-hand side represents the exponentially small contribution, where
the coefficients $A_j$ and $B_j$ are given by
\bee\label{e32}
A_j=\frac{ (1-a)_j(b-a)_j}{j!},\qquad B_j=\sum_{k=0}^j (-2)^{k} (\fs)_k\,A_{j-k}\,G_{2k,j-k}
\ee
for $j\geq0$.
The first five even-order coefficients $G_{2k,j}\equiv 6^{-2k} {\hat G}_{2k,j}$ are\footnote{There was a misprint in the first term in ${\hat G}_{6,j}$ in \cite{P13}, which appeared as $-3226$ instead of $-3626$. This was pointed out by T. Pudlik \cite{TP}. The correct value was used in the numerical calculations described in \cite{P13}.}
\begin{eqnarray}
{\hat G}_{0,j}\!\!&=&\!\!\f{2}{3}-\gamma_j,\qquad {\hat G}_{2,j}=\f{1}{15}(46-225\gamma_j+270\gamma_j^2-90\gamma_j^3), \nonumber\\
{\hat G}_{4,j}\!\!&=&\!\!\f{1}{70}(230-3969\gamma_j+11340\gamma_j^2-11760\gamma_j^3+5040\gamma_j^4
-756\gamma_j^5),\nonumber\\
{\hat G}_{6,j}\!\!&=&\!\!\f{1}{350}(-3626-17781\gamma_j+183330\gamma_j^2-397530\gamma_j^3+370440\gamma_j^4
-170100\gamma_j^5\nonumber\\
&&\hspace{7cm}+37800\gamma_j^6-3240\gamma_j^7),\nonumber\\
{\hat G}_{8,j}\!\!&=&\!\!\f{1}{231000}(-4032746+43924815\gamma_j+88280280\gamma_j^2-743046480\gamma_j^3\nonumber\\
&&+1353607200\gamma_j^4-1160830440\gamma_j^5+541870560\gamma_j^6
-141134400\gamma_j^7\nonumber\\
&&\hspace{6cm}+19245600\gamma_j^8-1069200\gamma_j^9)\label{e35}
\end{eqnarray}
with
\bee\label{e36}
\gamma_j:=\alpha-j\qquad (1\leq j\leq M-1).
\ee
The procedure for the generation of the coefficients $G_{2k,j}$ is given in the appendix.
From this it is evident that the coefficients $B_j$ not only depend on $a$ and $b$ but also on $\alpha$ in (\ref{e33}), which in turn depends on the particular value of the variable $x$ under consideration.

We can now apply the expansion in (\ref{e31}) to the representation of the solutions $w_1(a,x)$ and $w_2(a,x)$ in terms of Kummer functions in (\ref{e12}). With $a\to \fs a+\f{1}{4}$, $b=\fs$, the coefficients $A_j$ are
\[A_j=\frac{(\f{1}{4}-\fs a)_j (\f{3}{4}-\fs a)_j}{j!}=\frac{(\fs-a)_{2j}}{2^{2j} j!},\]
where we have employed the identity $(\alpha)_{2k}=2^{2k}(\fs\alpha)_k (\fs \alpha+\fs)_k$.
Then we obtain:
\newtheorem{theorem}{Theorem}
\begin{theorem}$\!\!\!.$\  The following expansions hold as $x\to+\infty$:
\bee\label{e37a}
w_1(a,x)-\frac{2^{\frac{1}{4}+\fr a}\sqrt{\pi}}{\g(\f{1}{4}-\fs a)}\,x^{-a-\fr} \sum_{k=0}^{m_o-1}\frac{(\fs+a)_{2k}}{k! (2x^2)^k}= E_1(a,x),\ee
\bee\label{e38a}
w_2(a,x)-\frac{2^{\fr a-\frac{1}{4}}\sqrt{\pi}}{\g(\f{3}{4}-\fs a)}\,x^{-a-\fr} \sum_{k=0}^{m_o-1}\frac{(\fs+a)_{2k}}{k! (2x^2)^k}= E_2(a,x),\ee
where, for positive integer $M$, the exponentially small expansions $E_1(a,x)$ and $E_2(a,x)$ are given by
\[E_1(a,x)=\frac{2^{\frac{1}{4}-\fr a}\sqrt{\pi}}{\g(\f{1}{4}+\fs a)}\,x^{a-\fr} e^{-x^2/2}\bl\{\cos \pi\vartheta \sum_{j=0}^{M-1} \frac{(-)^j (\fs-a)_{2j}}{j! (2x^2)^j}\hspace{4cm} \]
\bee\label{e37b}
\hspace{6cm}-\frac{2\sin \pi\vartheta}{\sqrt{\pi}}  \sum_{j=0}^{M-1}\frac{(-)^j 2^j B_j}{x^{2j+1}}+O(x^{-2M})\br\},\ee
\[E_2(a,x)=\frac{2^{-\fr a-\frac{1}{4}}\sqrt{\pi}}{\g(\f{3}{4}+\fs a)}\,x^{a-\fr}e^{-x^2/2}\bl\{\sin \pi\vartheta \sum_{j=0}^{M-1} \frac{(-)^j (\fs-a)_{2j}}{j! (2x^2)^j} \hspace{4cm}\]
\bee\label{e38b}
\hspace{6cm}+\frac{2\cos \pi\vartheta}{\sqrt{\pi}}  \sum_{j=0}^{M-1}\frac{(-)^j 2^j B_j}{x^{2j+1}}+O(x^{-2M})\br\}.\ee 
The parameter $\vartheta=\fs a-\f{1}{4}$, $m_o=\fs x^2-a+\alpha$, $|\alpha|<1$, is the optimal truncation index of the dominant algebraic expansion and the coefficients $B_j$ are given in (\ref{e32}) and (\ref{e35}) with the quantity $\gamma_j=\alpha-j$.
\end{theorem}
 
Comparison of the above expansions with those in (\ref{e24a}) and (\ref{e24b}) reveals two important differences. First, the factors $2\cos \pi\vartheta$ and $2\sin \pi\vartheta$ present in (\ref{e24a}) and (\ref{e24b}) have become $\cos \pi\vartheta$ and $\sin \pi\vartheta$, respectively. And secondly, with the dominant algebraic expansion optimally truncated, there appears an additional contribution to the exponentially small component given by the series involving the coefficients $B_j$.
\vspace{0.6cm}

\begin{center}
{\bf 3.\  Numerical examples and concluding remarks}
\end{center}
\setcounter{section}{3}
\setcounter{equation}{0}
\renewcommand{\theequation}{\arabic{section}.\arabic{equation}}
In this section we present some numerical examples to demonstrate the accuracy of the expansions in 
Theorem 1. As has already been noted, the coefficients $B_j$ depend on the parameter $a$ and also on $\alpha$ (see the definition of $\gamma_j$ in (\ref{e36})), which appears in the value of the optimal truncation index $m_o$ in (\ref{e33}). The value of $\alpha$ clearly is a function of the particular value of $x$ being considered.
In Table 1 we show values of the coefficients $A_j$ and $B_j$ for $0\leq j\leq 5$ and two different values of $a$ and $\alpha$. 
\begin{table}[h]
\caption{\footnotesize{Values of the coefficients $A_j$ and $B_j$ for $0\leq j\leq5$. }}
\begin{center}
\begin{tabular}{l|cc|cc}
\hline
&&&&\\[-0.30cm]
\mcol{1}{c|}{} & \mcol{2}{c|}{$a=1/4,\ \alpha=1/4$} & \mcol{2}{c}{$a=5/4,\ \alpha=0$}\\
\mcol{1}{c|}{$j$} & \mcol{1}{c}{$A_j$} & \mcol{1}{c|}{$B_j$} & \mcol{1}{c}{$A_j$} & \mcol{1}{c}{$B_j$}\\
[.1cm]\hline
&&&&\\[-0.30cm]
0 & 1.0000000000 & 0.4166666667 & $+1.0000000000$ & $+0.6666666667$\\
1 & 0.0781250000 & 0.1010127315 & $-0.0468750000$ & $-0.1633101852$\\
2 & 0.0714111328 & 0.1068229877 & $-0.0164794922$ & $+0.0184348132$\\
3 & 0.1327800751 & 0.2659511653 & $-0.0189685822$ & $-0.0474528804$\\
4 & 0.3760373220 & 0.8932217131 & $-0.0389004126$ & $-0.0734894988$\\
5 & 1.4348174067 & 3.8427298888 & $-0.1163365465$ & $-0.2841972836$\\
[.2cm]\hline
\end{tabular}
\end{center}
\end{table}

We define the quantities
\[{\cal W}_1(a,x):=w_1(a,x)-\frac{2^{\frac{1}{4}+\fr a}\sqrt{\pi}}{\g(\f{1}{4}-\fs a)}\,x^{-a-\fr} \sum_{k=0}^{m_o-1}\frac{(\fs+a)_{2k}}{k! (2x^2)^k},\]
\[{\cal W}_2(a,x):=w_2(a,x)-\frac{2^{\fr a-\frac{1}{4}}\sqrt{\pi}}{\g(\f{3}{4}-\fs a)}\,x^{-a-\fr} \sum_{k=0}^{m_o-1}\frac{(\fs+a)_{2k}}{k! (2x^2)^k},\]
which correspond to subtraction of the dominant, optimally truncated algebraic expansion from the solutions $w_1(a,x)$ and $w_2(a,x)$. In Table 2 we show the values\footnote{In Tables 2 and 3 we write the values as $x(y)$ instead of $x\times 10^y$.} of the exponentially small expansions $E_1(a,x)$ and $E_2(a,x)$ for different truncation index $M$ together with the computed values of ${\cal W}_1(a,x)$ and ${\cal W}_2(a,x)$. Finally, in Table 3 the values of ${\cal W}_r(ax)$ and $E_r(a,x)$ are presented for different values of $x$ using the truncation index $M=6$. These results indicate the validity of the expansions in Theorem 1.
\begin{table}[h]
\caption{\footnotesize{Values of $E_1(a,x)$ and $E_2(a,x)$ for different truncation index $M$ when $x=6$.
The values of ${\cal W}_r(a,x)$ ($r=1, 2$) are shown for comparison in the last row.}}
\begin{center}
\begin{tabular}{c|cc|cc}
\hline
&&&&\\[-0.30cm]
\mcol{1}{c|}{} & \mcol{2}{c|}{$a=1/4,\ \ \alpha=1/4$} & \mcol{2}{c}{$a=1,\ \ \alpha=0$}\\
\mcol{1}{c|}{$M$} & \mcol{1}{c}{$E_1(a,x)$} & \mcol{1}{c|}{$E_2(a,x)$} & \mcol{1}{c}{$E_1(a,x)$} & \mcol{1}{c}{$E_2(a,x)$}\\
[.1cm]\hline
&&&&\\[-0.30cm]
1 & $7.5687562(-9)$ & $-3.7792328(-9)$ & $2.8061782(-8)$ & $3.4517454(-8)$\\
2 & $7.5337338(-9)$ & $-3.7789129(-9)$ & $2.8109713(-8)$ & $3.4684633(-8)$\\
3 & $7.5355378(-9)$ & $-3.7792432(-9)$ & $2.8106765(-8)$ & $3.4681904(-8)$\\
4 & $7.5353448(-9)$ & $-3.7792335(-9)$ & $2.8106892(-8)$ & $3.4682187(-8)$\\
5 & $7.5353760(-9)$ & $-3.7792322(-9)$ & $2.8106876(-8)$ & $3.4682152(-8)$\\
6 & $7.5353692(-9)$ & $-3.7792330(-9)$ & $2.8106879(-8)$ & $3.4682159(-8)$\\
[.15cm]\hline
&&&&\\[-0.30cm]
${\cal W}_r(a,x)$ & $7.5353706(-9)$ & $-3.7792328(-9)$ & $2.8106878(-8)$ & $3.4682157(-8)$\\
\end{tabular}
\end{center}
\end{table}
\begin{table}[h]
\caption{\footnotesize{Values of ${\cal W}_r(a,x)$ and $E_r(a,x)$ ($r=1, 2$) for different $x$ when $a=1$ and the truncation index $M=6$.}}
\begin{center}
\begin{tabular}{c|cc|cc}
\hline
&&&&\\[-0.30cm]
%\mcol{1}{c|}{} & \mcol{2}{c||}{$a=1/4,\ \ \alpha=-1/4$} & \mcol{2}{c}{$a=1,\ \ \alpha=0$}\\
\mcol{1}{c|}{$x$} & \mcol{1}{c}{${\cal W}_1(a,x)$} & \mcol{1}{c|}{$E_1(a,x)$} & \mcol{1}{c}{${\cal W}_2(a,x)$} & \mcol{1}{c}{$E_2(a,x)$}\\
[.1cm]\hline
&&&&\\[-0.30cm]
2 & $9.7647365111(-02)$ & $9.7{\bf 2}32594660(-02)$ & $2.2968994497(-01)$ & $2.{\bf 3}181240997(-01)$\\
3 & $1.5656185046(-02)$ & $1.56561{\bf 9}3695(-02)$ & $1.7075649565(-02)$ & $1.707{\bf 7}772223(-02)$\\
4 & $4.6890418500(-04)$ & $4.68904{\bf 5}6968(-04)$ & $6.6631593463(-04)$ & $6.6631{\bf 6}19766(-04)$\\
5 & $6.9251877004(-06)$ & $6.92518{\bf 8}3019(-06)$ & $7.1483664282(-06)$ & $7.14836{\bf 8}7955(-06)$\\
6 & $2.8106878174(-08)$ & $2.8106878{\bf 5}86(-08)$ & $3.4682157382(-08)$ & $3.468215{\bf 8}600(-08)$\\
8 & $2.7943166845(-14)$ & $2.79431668{\bf 6}2(-14)$ & $3.2300734782(-14)$ & $3.2300734782(-14)$\\
[.15cm]\hline
\end{tabular}
\end{center}
\end{table}

\vspace{0.6cm}
\newpage
\begin{center}
{\bf Appendix: \ The generation of the coefficients $G_{2k,j}$}
\end{center}
\setcounter{section}{1}
\setcounter{equation}{0}
\renewcommand{\theequation}{\Alph{section}.\arabic{equation}}
The coefficients $G_{2k,j}$ are generated from the expansion (see \cite[p.~1468]{O})
\bee\label{a1}
\frac{t^{\gamma_j-1}}{1-t}\,\frac{dt}{dw}+\frac{1}{w}=\sum_{k=0}^\infty G_{2k,j}w^k,\qquad \fs w^2=t-\log\,t-1,\ee
where $\gamma_j$ is defined in (\ref{e36}). The branch of $w(t)$ is chosen such that $w\sim t-1$ as $t\to 1$. Reversion if the $w-t$ mapping yields
\[t=1+w+\f{1}{3}w^2+\f{1}{36}w^3-\f{1}{270}w^4+\f{1}{4320}w^5+\cdots\ .\]
Substitution of this last relation into (\ref{a1}) then enables with the help of {\it Mathemtica} the calculation of the even coefficients $G_{2k,j}$.

\end{document}